\documentclass[12pt,draft]{amsart}
\usepackage{amsmath,amsthm,latexsym,amscd,amsbsy,amssymb,fleqn,leqno}
\chardef\bslash=`\\ % p. 424, TeXbook

\makeatletter
\def\verbatim{\interlinepenalty\@M \@verbatim
  \leftskip\@totalleftmargin\advance\leftskip2pc
  \frenchspacing\@vobeyspaces \@xverbatim}
\makeatother
\makeatletter
  \def\dgt@k{\dg@DX=-3 \dg@DY=2 \dg@SIZE=3} 
\makeatother

\makeatletter
  \def\dgt@kk{\dg@DX=3 \dg@DY=-1 \dg@SIZE=3}%
\makeatother

\theoremstyle{plain}
\newtheorem{thm}{Theorem}[section]
\newtheorem{cor}[thm]{Corollary}
\newtheorem{lem}[thm]{Lemma}
\newtheorem{pro}[thm]{Proposition}

\theoremstyle{definition}

\numberwithin{equation}{section}

\newcounter{rmnum}

\def\symbolnote#1#2{\let\thefootn=\thefootnote%
\renewcommand{\thefootnote}{\fnsymbol{footnote}}%
\footnotemark[#1]%
\footnotetext[#1]{#2}%
\let\thefootnote=\thefootn
}

\newfont{\bbb}{msbm10 scaled \magstep1}
\newfont{\bbc}{msbm8 scaled \magstep0}

\newcommand{\R}{\mbox{\bbb R}}

\newcommand{\N}{\mbox{\bbb N}}

\newcommand{\sphere}{\mbox{\bbb S}}
\newcommand{\uin}{\mbox{\bbb I}}
\newcommand{\e}{\mbox{\rm e-dim}}

\begin{document}

%%%%%%% Begin Topmatter %%%%%%%%%%

\title[On dimensionally restricted maps]{On dimensionally restricted maps}
\author{H. Murat Tuncali}
\address{Department of Mathematics,
Nipissing University\-wan,
100 College Drive, P.O. Box 5002, North Bay, ON, P1B 8L7, Canada}
\email{muratt@unipissing.ca}
\thanks{The first author was partially supported by NSERC grant.}

\author{Vesko Valov}
\address{Department of Mathematics, Nipissing University,
100 College Drive, P.O. Box 5200, North Bay, ON, P1B 8L7, Canada}
\email{veskov@unipissing.ca}
\thanks{The second author was partially supported by Nipissing University Research Council Grant.}

\keywords{finite-dimensional maps, extensional dimension, $C$-space} 
\subjclass{Primary: 54F45; Secondary: 55M10, 54C65.}
 
%%%%%%% End topmatter %%%%%%%%%

\begin{abstract}{Let $f\colon X\to Y$ be a closed $n$-dimensional surjective map of metrizable spaces. It is shown that if $Y$ is a $C$-space, then: (1) the set of all maps $g\colon X\to\uin^n$ with $\dim (f\times g)=0$ is uniformly dense in $C(X,\uin^n)$; (2) for every $0\leq k\leq n-1$ there exists an $F_{\sigma}$-subset $A_k$ of $X$ such that $\dim A_k\leq k$ and the restriction $f|(X\backslash A_k)$ is $(n-k-1)$-dimensional. These are extensions of theorems by Pasynkov and Torunczyk, respectively, obtained for finite-dimensional spaces. A generalization of a result due to Dranishnikov and Uspenskij about extensional dimension is also established.}
\end{abstract}

\maketitle
\markboth{H. M.~Tuncali and V.~Valov}{On dimensionally restricted maps}

%%%%%%%%%%%%%%%%%%%%%%%%%%%%%%%%%%%%%%%%%%%%%%%%%%%%%%%%%%%%

\section{Introduction}

All spaces are assumed to be completely regular and all maps continuous.
This paper concerns with the following two results. The first one was proved by Pasynkov  \cite{bp:96} (see \cite{bp:98} for non-compact versions) and the second one by Torunczyk \cite{ht:85}:

\begin{thm}$($Pasynkov$)$.
Let $f\colon X\to Y$ be an $n$-dimensional map with $X$ and $Y$ being finite-dimensional compact metric spaces. Then there exists $g\colon X\to \uin^n$ such that 
$f\times g\colon X\to Y\times\uin^n$ is $0$-dimensional. Moreover, the set of all such $g$ is dense and $G_{\delta}$ in $C(X,\uin^n)$ with respect to uniform convergence topology. 
\end{thm}

\begin{thm}$($Torunczyk$)$.
Let $f\colon X\to Y$ be a $\sigma$-closed map of separable metric spaces with $\dim f=n$ and $\dim Y<\infty$. Then for each $0\leq k\leq n-1$ there exists an $F_{\sigma}$-subset $A_k$ of $X$ such that $\dim A_k\leq k$ and the restriction $f|(X\backslash A_k)$ is $(n-k-1)$-dimensional.
\end{thm}

The above two theorems are equivalent in the realm of compact spaces (see \cite{ml:96} and \cite{ys:95}). However, the problem whether these two theorems hold without any dimensional restrictions on $Y$ is still open. Sternfeld and Levin made a significant progress in solving this problem.
In 1995, Sternfeld \cite{ys:95} proved that if $f\colon X\to Y$ is an $n$-dimensional map between compact metric spaces, then $\dim (f\times g)\leq 1$ for almost all $g\in C(X,\uin^n)$; equivalently, there exists a $\sigma$-compact $(n-1)$-dimensional subset $A$ of $X$ such that $\dim(f|(X\backslash A))\leq 1$.
Levin \cite{ml:96} improved Sternfeld's result
showing that $\dim (f\times g)\leq 0$ for almost all maps $g\in C(X,\uin^{n+1})$ which is equivalent to the existence of a $\sigma$-compact $n$-dimensional set $A\subset X$ with $\dim(f|(X\backslash A))\leq 0$. 

In the present paper we generalize Theorem 1.1 and Theorem 1.2 to arbitrary metrizable spaces by replacing the finite dimensionality of $Y$ with the less restrictive condition to be a $C$-space. Recall that a space $X$ is a $C$-space \cite{ag:78} if for any sequence
$\{\omega_n:n\in\N\}$ of open covers of $X$ there exists a sequence $\{\gamma_n:n\in\N\}$ of open disjoint families in $X$ such that each $\gamma_n$ refines $\omega_n$ and $\cup\{\gamma_n:n\in\N\}$ covers $X$. 
$C$-space property was introduced by Haver \cite{wh:74} for compact metric spaces and Addis and Gresham \cite{ag:78} extended Haver's definition for more general spaces.  
All countable-dimensional metrizable spaces (spaces which are countable union of finite-dimensional subsets), in particular all finite-dimensional ones, form a proper subclass of the class of $C$-spaces because there exists a metric $C$-compactum which is not countable-dimensional \cite{p:81}. 

Here is the generalized version of Theorem 1.1.

\begin{thm}
Let $f\colon X\to Y$ be a closed map of metric spaces with $\dim f=n$ and $Y$ a $C$-space. Then all maps $g\colon X\to \uin^n$ such that $\dim(f\times g)=0$ form a dense subset of $C(X,\uin^n)$ with respect to the uniform convergence topology. Moreover, if $f$ is $\sigma$-perfect, then this set is dense and $G_{\delta}$ in $C(X,\uin^n)$ with respect to the source limitation topology.
\end{thm} 

Theorem 1.3 answers positively Pasynkov's question in \cite{bp:96} whether Theorem 1.1 is true for countable-dimensional spaces. 

For any map $f\colon X\to Y$ 
$\dim f=\sup\{\dim f^{-1}(y):y\in Y\}$ is the dimension of $f$. We say that 
a surjective map $f\colon X\to Y$ is called $\sigma$-closed (resp., $\sigma$-perfect) if $X$ is the union of countably many closed sets $X_i$ such that 
each restriction $f|X_i\colon X_i\to f(X_i)$ is a closed (resp., perfect) map and all $f(X_i)$ are closed in $Y$. 

Using Theorem 1.3 we prove the following generalization of Theorem 1.2:

\begin{thm}
Let $f\colon X\to Y$ be a $\sigma$-closed map of metric spaces with $\dim f=n$ and $Y$ a $C$-space. Then for each $0\leq k\leq n-1$ there exists an $F_{\sigma}$-subset $A_k$ of $X$ such that $\dim A_k\leq k$ and the restriction $f|(X\backslash A_k)$ is $(n-k-1)$-dimensional.
\end{thm}

A few words about this note. In Section 2 we give a characterization of finite-dimensional  proper maps (see Theorem 2.2) which is the main tool in the proof of Theorem 1.3. 
The proof of Theorem 2.2 is based on a selection theorem established by V. Gutev and the second author \cite[Theorem 1.1]{gv:99}.
Sections 3 and 4 are devoted to the proof of Theorem 1.3 and Theorem 1.4, respectively.
In the last Section 5 we provide applications of the main results. One of them is a generalization of a result by Dranishnikov and Uspenskij \cite{du:97} concerning maps which lower extensional dimension, another one is a parametric version of the Bogatyi representation theorem of $n$-dimensional metrizable spaces \cite{sb:99}. Some results in the spirit of Pasynkov's recent paper \cite{bp:98} are also obtained.

%%%%%%%%%%%%%%%%%%%%%%%%%%%%%%%%%

%%%%%%%%%%%%%%%%%%%%%%%%%%%%%%
%%%%%%%%%%%%%%%%%%%%%%%%%%%%%%

\section{Finite-dimensional maps}
In this section we provide a characterization of $n$-dimensional perfect maps onto paracompact $C$-spaces, see Theorem 2.2 below. 
 
For any spaces $M$ and $K$ by $C(K,M)$ we denote the set of all continuous maps from
$K$ into $M$. If $(M,d)$ is a metric space and $K$ is any space, then the source limitation topology on $C(K,M)$ is defined in the following way:
a subset $U\subset C(K,M)$ is open in $C(K,M)$ with respect to
the source limitation topology provided for every $g\in U$ there exists 
a continuous function $\alpha\colon K\to (0,\infty)$ such that $\overline{B}(g,\alpha)\subset U$. Here, $\overline{B}(g,\alpha)$ denotes the set  
$\{h\in C(K,M):d(g(x),h(x))\leq\alpha (x)\hbox{}~~\mbox{for each 
$x\in K$}\}$. 

The source limitation topology is also known as the fine topology and $C(K,M)$ with this topology has Baire property provided $(M,d)$ is a complete metric space \cite{jm:75}.  
We also need the following fact: if $K$ is paracompact and $F\subset K$ closed, then 
the restriction map $p_F\colon C(K,M)\to C(F,M)$, $p_F(g)=g|F$, is continuous when both
$C(K,M)$ and $C(F,M)$ are equipped with the source limitation topology; moreover $p_F$ is open and surjective provided $M$ is a closed convex subset of a Banach space and $d$ is the metric on $M$ generated by the norm. Finally, when $K$ and $M$ are metrizable, the source limitation topology on $C(K,M)$ doesn't depend on the concrete metric on $M$. 

Let  $\omega$ be an open cover of the space $M$ and $m\in\N\cup\{0\}$. A family $\gamma$ of subsets of $M$ is said to be $(m,\omega)$-discrete 
in $M$ if
$ord(\gamma)\leq m+1$ (i.e., every point of $M$ belongs to at most $m+1$ elements of $\gamma$) and $\gamma$ refines $\omega$; a subset of $M$ which can be represented as the union of open $(m,\omega)$-discrete family in $M$ is called $(m,\omega)$-discrete; a map
$g\colon M\to Z$ is $(m,\omega)$-discrete if every $z\in g(M)$ has a neighborhood $V_z$ in $Z$ such that $g^{-1}(V_z)$ is $(m,\omega)$-discrete in $M$. 

We also agree to denote by $cov(M)$ the family of all open covers of $M$. In case $(M,d)$ is a metric space,
$B_{\epsilon}(x)$ (resp., $\overline{B}_{\epsilon}(x)$) stands for the
open (resp., closed) ball in $(M,d)$ with center $x$ and radius
$\epsilon$.

\begin{lem}
If $\omega\in cov(M)$ and $K\subset M$ compact, then every functionally open and $(m,\omega)$-discrete subset of $K$ can be extended to an $(m,\omega)$-discrete subset of $M$.
\end{lem}

\begin{proof}
Let $U\subset K$ be functionally open and $(m,\omega)$-discrete in $K$ and  $\gamma =\{U_s:s\in A\}$ an open $(m,\omega)$-discrete family in $K$ whose union is $U$. Since $U$ is paracompact (as functionally open in $K$), we can suppose that $\gamma$ is locally finite and there exists a partition of unity $\{f_s:s\in A\}$ in $U$ such that $U_s=f_s^{-1}((0,1])$ for each $s\in A$. Denote by $\mathcal N$ the nerve of $\gamma$ with the Whitehead topology and define the map $f\colon U\to {\mathcal N}$ by $f(x)=\sum_{s\in A}f_s(x)s$. Observe that $\mathcal N$ is at most $m$-dimensional because $ord(\gamma)\leq m+1$.
Let $W$ be a functionally open subset of $\beta M$ with $W\cap K=U$. Then, by \cite{jd:99}, there exists an open set $V\subset W$ containing $U$ and an extension $g\colon V\to\mathcal N$ of $f$. The map $g$ generates maps $g_s\colon V\to [0,1]$ such that each $g_s$ extends $f_s$. We finally choose $G_s\in\omega$ with $U_s\subset G_s$, $s\in A$, and define $V_s=G_s\cap g_s^{-1}((0,1])$. Then  
the family $\{V_s:s\in A\}$ is $(m,\omega)$-discrete in $M$ and $\bigcup_{s\in A}V_s$ is the required $(m,\omega)$-discrete extension of $U$.
\end{proof} 
 
Throughout the paper $\uin^k$ denotes the $k$-dimensional cube equipped with the Euclidean metric $d_k$, and $D_k$ denotes the uniform convergence metric on $C(X,\uin^k)$ generated by $d_k$. 
If $f\colon X\to Y$, we denote by $C(X,Y\times\uin^k,f)$
the set of all maps $h\colon X\to Y\times\uin^k$ such that $\pi_Y\circ h=f$, where
$\pi_Y\colon Y\times\uin^k\to Y$ is the projection. For any $\omega\in cov(X)$ and closed $K\subset X$, $C_{(m,\omega)}(X|K,Y\times\uin^k,f)$ stands for the set of all $h\in C(X,Y\times\uin^k,f)$
with $h|K$ being $(m,\omega)$-discrete (as a map from $K$ into $Y\times\uin^k$) and 
$C_{(m,\omega)}(X|K,\uin^k)$ consists of all $g\in C(X,\uin^k)$ such that 
$f\times g\in C_{(m,\omega)}(X|K,Y\times\uin^k,f)$. In case $K=X$ we simply write 
$C_{(m,\omega)}(X,Y\times\uin^k,f)$ (resp., $C_{(m,\omega)}(X,\uin^k)$) instead of
$C_{(m,\omega)}(X|X,Y\times\uin^k,f)$ (resp., $C_{(m,\omega)}(X|X,\uin^k)$). 

Now we can establish the following characterization of $n$-dimensional perfect maps:

\begin{thm}
Let $f\colon X\to Y$ be a perfect surjection between paracompact spaces with $Y$ being a $C$-space. Then $\dim f\leq n$ if and only if for any $\omega\in cov(X)$ and $0\leq k\leq n$ the set
$C_{(n-k,\omega)}(X,\uin^k)$ is open and dense in $C(X,\uin^k)$ with respect to the source limitation topology.
\end{thm}

The proof of sufficiency follows from the following observation: if the set $C_{(0,\omega)}(X,\uin^n)$ is not empty for all $\omega\in cov(X)$, then every open cover of $f^{-1}(y)$, $y\in Y$, admits an open refinement of order $\leq n+1$, i.e. $\dim f^{-1}(y)\leq n$. 
Indeed, let $\gamma$ be a family of open subsets of $X$ covering $f^{-1}(y)$. Then $\omega =\gamma\cup\{X\backslash f^{-1}(y)\}\in cov(X)$, so there exists 
$g\in C_{(0,\omega)}(X,\uin^n)$. Obviously, $g|f^{-1}(y)$ is $(0,\omega)$-discrete. Hence, every $z\in H=g(f^{-1}(y))$ has a neighborhood $G_z$ in $\uin^n$ with $g^{-1}(G_z)\cap f^{-1}(y)$ being the union of a disjoint and open in $f^{-1}(y)$ family $\mu_z$ which refines $\omega$. Take finitely many $z(i)\in H$, $i=1,2,..,p$, such that $\lambda=\{G_{z(i)}:i=1,2,..,p\}$ covers $H$. Since $\dim H\leq n$, we can suppose that $ord(\lambda)\leq n+1$. Then
$\mu=\cup\{\mu_{z(i)}:i=1,2,..,p\}$ is an open cover of $f^{-1}(y)$ refining $\gamma$ and $ord(\mu)\leq n+1$.
  
To prove necessity we need few lemmas, the proof will be completed by Lemma 2.9. 
In all these lemmas we suppose that $X$ and $Y$ are given paracompact spaces and
$f\colon X\to Y$ a perfect surjective map with $\dim f\leq n$, where $n\in\N$. We also fix $\omega\in cov(X)$, an integer $k$ such that $0\leq k\leq n$
and arbitrary $m\in\N\cup\{0\}$. 

\begin{lem}
Let $g\in C_{(m,\omega)}(X|f^{-1}(y),\uin^k)$ for some $y\in Y$. Then there exists a neighborhood $U$ of $y$ in $Y$ such that the restriction $g|f^{-1}(U)$ is $(m,\omega)$-discrete.
\end{lem}

\begin{proof}
Obviously, $g\in C_{(m,\omega)}(X|f^{-1}(y),\uin^k)$ implies that 
$g|f^{-1}(y)$ is an $(m,\omega)$-discrete map. Hence, for every $x\in f^{-1}(y)$ there exists an open neighborhood $V_{g(x)}$ of $g(x)$ in $\uin^k$ such that $g^{-1}(V_{g(x)})\cap f^{-1}(y)$ is an $(m,\omega)$-discrete set in $f^{-1}(y)$. Since $V_{g(x)}$ is functionally open in 
$\uin^k$, so is $g^{-1}(V_{g(x)})\cap f^{-1}(y)$ in $f^{-1}(y)$. Then, by Lemma 2.1, there is
an $(m,\omega)$-discrete subset $W_x$ in $X$ extending $g^{-1}(V_{g(x)})\cap f^{-1}(y)$.
Therefore, for every $x\in f^{-1}(y)$ we have
$(f\times g)^{-1}(f(x),g(x))=f^{-1}(y)\cap g^{-1}(g(x))\subset W_x$ and, since $f\times g$ is a closed map, there exists an open neighborhood $H_x=U_{y}^x\times G_x$ of $(y,g(x))$ in $Y\times\uin^k$ with $S_x=(f\times g)^{-1}(H_x)\subset W_x$. Next, choose finitely many points $x(i)\in f^{-1}(y)$, $i=1,2,..p$, 
such that $f^{-1}(y)\subset \bigcup_{i=1}^{i=p}S_{x(i)}$. Using that $f$ is a closed map we can find a neighborhood $U_y$ of $y$ in $Y$ such that $U_y\subset\bigcap_{i=1}^{i=p}U_{y}^{x(i)}$ and $f^{-1}(U_y)\subset\bigcup_{i=1}^{i=p}S_{x(i)}$. Let show that $g|f^{-1}(U_y)$ is $(m,\omega)$-discrete. Indeed, if $z\in f^{-1}(U_y)$, then $z\in S_{x(j)}$ for some $j$ and $g(z)\in G_{x(j)}$ because $S_{x(j)}=f^{-1}(U_{y}^{x(j)})\cap g^{-1}(G_{x(j)})$. Consequently, 
$f^{-1}(U_{y})\cap g^{-1}(G_{x(j)})\subset S_{x(j)}\subset W_{x(j)}$. Therefore, $G_{x(j)}$ is a neighborhood of $g(z)$ such that $f^{-1}(U_{y})\cap g^{-1}(G_{x(j)})$ is $(m,\omega)$-discrete in $f^{-1}(U_y)$ as a subset of the $(m,\omega)$-discrete set $W_{x(j)}$ in $X$.
\end{proof}

\begin{cor}
If $g\in C_{(m,\omega)}(X|f^{-1}(y),\uin^k)$ for every $y\in Y$, then we have
$g\in C_{(m,\omega)}(X,\uin^k)$.
\end{cor}

\begin{proof}
We need to show that $f\times g$ is $(m,\omega)$-discrete, i.e. for any $x\in X$ there exist
neighborhoods $U_y$ of $y=f(x)$ in $Y$ and $G_x$ of $g(x)$ in $\uin^k$ such that 
$f^{-1}(U_y)\cap g^{-1}(G_x)$ is $(m,\omega)$-discrete in $X$. And this is really true, 
by Lemma 2.3, there exists a neighborhood $U_y$ of $y$ in $Y$ such that $g|f^{-1}(U_y)$ is $(m,\omega)$-discrete. Therefore, we can find a neighborhood $G_{x}$ of $g(x)$ in $\uin^k$ with  
$f^{-1}(U_{y})\cap g^{-1}(G_{x})$ being $(m,\omega)$-discrete in $f^{-1}(U_y)$. Consequently, $f^{-1}(U_{y})\cap g^{-1}(G_{x})$ is $(m,\omega)$-discrete in $X$.
\end{proof} 

\begin{lem}
The set 
$C_{(m,\omega)}(X|K,\uin^k)$ is open in $C(X,\uin^k)$ with respect to the source limitation topology for any closed $K\subset X$.  
\end{lem}

\begin{proof}
Let $g_0\in C_{(m,\omega)}(X|K,\uin^k)$. We are going to find $\alpha\in C(X,(0,\infty))$ with $\overline{B}(g_0,\alpha)\subset C_{(m,\omega)}(X|K,\uin^k)$. Since each restriction $g_0|(f^{-1}(y)\cap K)$, $y\in H=f(K)$, is $(m,\omega)$-discrete, by Lemma 2.3, for every $y\in H$ there exists a neighborhood $U_y$  of $y$ in $Y$ such that $g_0|(f^{-1}(U_y)\cap K)$ is $(m,\omega)$-discrete. Then $\omega_1=\{U_y:y\in H\}\cup\{Y\backslash H\}$ is an open cover of $Y$. Using that $Y$ is paracompact, we can find a metric space $(M,d)$, a surjection $p\colon Y\to M$ and $\mu\in cov(M)$ such that $p^{-1}(\mu)$ refines $\omega_1$.
Hence, every $z\in p(H)$ has a neighborhood $W_z$ in $M$ such that $g_0|(p\circ f)^{-1}(W_z)\cap K$ is $(m,\omega)$-discrete. The last condition implies that $h_0|K$ is $(m,\omega)$-discrete, where 
$h_0=(p\circ f)\times g_0$. Now we need the following:

\medskip
{\em Claim.} {\em There exists an open family $\gamma$ in $M\times\uin^k$ covering $h_0(K)$ such that
every $g\in C(X,\uin^k)$ belongs to $C_{(m,\omega)}(X|K,\uin^k)$ provided $h|K$ is $\gamma$-close to $h_0|K$, where $h=(p\circ f)\times g$.}

\medskip 

Proof of the claim. Since $h_0|K$ is $(m,\omega)$-discrete, every $t\in h_0(K)$ has an open neighborhood $V_t$ in $M\times\uin^k$ such that 
$h_{0}^{-1}(V_t)\cap K$ is $(m,\omega)$-discrete in $K$. 
Then $\nu =\{V_t:t\in h_0(K)\}$ forms an open cover of $h_0(K)$. Take $\gamma$ to be a locally finite open cover of $V=\cup\nu$ such that $\{St(W,\gamma):W\in\gamma\}$ refines $\nu$. 
Let $h|K$ be a $\gamma$-close map to $h_0|K$, where $h=(p\circ f)\times g$ with $g\in C(X,\uin^k)$. If $W\in\gamma$, then
$h_0(h^{-1}(W)\cap K)\subset St(W,\gamma)$. But $St(W,\gamma)$ is contained in $V_t$ for some $t\in h_0(K)$. Consequently, $h^{-1}(W)\cap K\subset h_0^{-1}(V_t)\cap K$. 
The last inclusion implies that $h^{-1}(W)\cap K$ is $(m,\omega)$-discrete in $K$ because $h_0^{-1}(V_t)\cap K$ is. Therefore, $h|K$ is $(m,\omega)$-discrete. To finish the proof of the claim observe that $h|K$ being $(m,\omega)$-discrete yields $(f\times g)|K$ is $(m,\omega)$-discrete too, i.e. $g\in C_{(m,\omega)}(X|K,\uin^k)$.  

\smallskip
We continue with the proof of Lemma 2.5. Let $\rho$ be the metric on $M\times\uin^k$ defined by $\rho (t_1,t_2)=d(z_1,z_2)+d_k(w_1,w_2)$, where $t_i=(z_i,w_i)$, $i=1,2$. 
Let $\alpha_1 \colon K\to (0,\infty)$ be the function 
$\alpha_1 (x)=2^{-1}\sup\{\rho (h_0(x),V\backslash W):W\in\gamma\}$. Since $h_0(K)\subset V$ and $\gamma$ is a locally finite open cover of $V$, $\alpha_1$ is continuous. Moreover, if 
$h=(p\circ f)\times g$ with $g\in C(X,\uin^k)$ and
$\rho(h_0(x),h(x))\leq\alpha_1(x)$ for every $x\in K$, then $h|K$ is $\gamma$-close to 
$h_0|K$. According to the claim, the last relation yields that 
$g\in C_{(m,\omega)}(X|K,\uin^k)$. We finally take a continuous extension $\alpha\colon X\to (0,\infty)$ of $\alpha_1$. Observe that $d_k(g_0(x),g(x))=\rho(h_0(x),h(x))$ for every $x\in X$. 
Therefore, $\overline{B}(g_0,\alpha)\subset C_{(m,\omega)}(X|K,\uin^k)$.
\end{proof}

\begin{lem}
If $C(X,\uin^k)$ is equipped with the uniform convergence topology, then the set-valued map $\displaystyle\psi_{(m,\omega)}\colon Y\to 2^{C(X,\uin^k)}$, defined by the formula
$\psi_{(m,\omega)}(y)=C(X,\uin^k)\backslash C_{(m,\omega)}(X|f^{-1}(y),\uin^k)$, has a closed graph.
\end{lem}

\begin{proof}
Let $G=\cup\{y\times\psi_{(m,\omega)}(y):y\in Y\}\subset Y\times C(X,\uin^k)$ be the graph of $\psi_{(m,\omega)}$ and $(y_0,g_0)\in (Y\times C(X,\uin^k))\backslash G$. We are going to show that $(y_0,g_0)$ has a neighborhood in $Y\times C(X,\uin^k)$ which doesn't meet $G$.
Since $(y_0,g_0)\not\in G$, $g_0\not\in\psi_{(m,\omega)}(y_0)$. Hence,  
$g_0\in C_{(m,\omega)}(X|f^{-1}(y_0),\uin^k)$ and, by Lemma 2.3, there exists a
neighborhood $U$ of $y_0$ in $Y$ with $g_0|f^{-1}(U)$ being $(m,\omega)$-discrete, in particular, $g_0\in C_{(m,\omega)}(X|f^{-1}(U),\uin^k)$. 
We can assume that $U\subset Y$ is closed, so is $f^{-1}(U)$ in $X$. Then, according to Lemma 2.5, $C_{(m,\omega)}(X|f^{-1}(U),\uin^k)$ is open in $C(X,\uin^k)$ with respect to the source limitation topology. Consequently, there exists a continuous positive function $\alpha$ on $X$ such that $\overline{B}(g_0,\alpha)$ is contained in $C_{(m,\omega)}(X|f^{-1}(U),\uin^k)$.
Since $f^{-1}(y_0)$ is compact, $2\delta =\min\{\alpha (x):x\in f^{-1}(y_0)\}>0$ and 
$H=\{x\in f^{-1}(U):\alpha (x)>\delta\}$ is a neighborhood of $f^{-1}(y_0)$. Therefore, there exists a closed neighborhood $V$ of $y_0$ in $Y$ with $f^{-1}(V)\subset H$ (we use again that $f$ is a closed map).
Let $B_{\delta}(g_0)$ be the open ball in $C(X,\uin^k)$ (with respect to the uniform metric $D_k$) with center $g_0$ and radius $\delta$. Since $W=V\times B_{\delta}(g_0)$ is a neighborhood of $(y_0,g_0)$ in $Y\times C(X,\uin^k)$, the following claim completes the proof.

\medskip
{\em Claim.} {\em $W\cap G=\emptyset$}

\medskip 
Suppose $(y,g)\in W\cap G$ for some $(y,g)\in Y\times C(X,\uin^k)$. Then, $y\in V$ and 

\medskip\noindent   
(1) \hbox{}~~~~~~$d_k(g(x),g_0(x))\leq\delta <\alpha (x)$ for every $x\in f^{-1}(V)$. 

\medskip\noindent
Let show that the existence of a function $g_1\in C(X,\uin^k)$ such that 

\medskip\noindent
(2) \hbox{}~~~~~~$g_1\in\overline{B}(g_0,\alpha)$ and $g_1|f^{-1}(V)=g|f^{-1}(V)$

\smallskip\noindent
provides a contradiction with the assumption $(y,g)\in W\cap G$. Indeed, $g_1\in\overline{B}(g_0,\alpha)$ yields 
$g_1\in C_{(m,\omega)}(X|f^{-1}(U),\uin^k)$ and, since $f^{-1}(y)\subset f^{-1}(U)$, we have 
$g_1\in C_{(m,\omega)}(X|f^{-1}(y),\uin^k)$. So, 
$g\in C_{(m,\omega)}(X|f^{-1}(y),\uin^k)$ because
$g_1|f^{-1}(y)=g|f^{-1}(y)$ (recall that $f^{-1}(y)\subset f^{-1}(V)$). 
On the other hand, $(y,g)\in G$ implies $g\in\psi_{(m,\omega)}(y)$, i.e. $g\not\in C_{(m,\omega)}(X|f^{-1}(y),\uin^k)$.

Therefore, the proof is reduced to find $g_1$ satisfying $(2)$. And this can be done by using the convex-valued selection theorem of Michael \cite{em:56}. Define the set-valued map $\Phi\colon X\to {\mathcal F}_c(\uin^k)$ by $\Phi (x)=g(x)$ if $x\in f^{-1}(V)$ and $\Phi (x)= \overline{B}_{\alpha(x)}(g_0(x))$ otherwise. Here, ${\mathcal F}_c(\uin^k)$ denotes the convex and closed subsets of $\uin^k$ and $\overline{B}_{\alpha(x)}(g_0(x))$ is the closed ball in $\uin^k$ with center $g_0(x)$ and radius $\alpha(x)$. By virtue of $(1)$, 
$g(x)\in\overline{B}_{\alpha(x)}(g_0(x))$ for all $x\in f^{-1}(V)$. The last condition, together with the definition of $\Phi$ outside $f^{-1}(V)$, implies that 
$\Phi$ is lower semi-continuous (i.e., 
$\{x\in X:\Phi(x)\cap O\neq\emptyset\}$ is open in $X$
for any open set $O\subset\uin^k$). Then, by mentioned above Michael's theorem, $\Phi$ admits a continuous selection $g_1$. Since $g_1(x)\in\Phi(x)$ for any $x\in X$, we have  
$g_1|f^{-1}(V)=g|f^{-1}(V)$ and $g_1\in\overline{B}(g_0,\alpha)$.  
\end{proof}

\begin{lem}
Let $K$ and $M$ be compact spaces such that $\dim K\leq n$ and $M$ metrizable.
Then for every $\gamma\in cov(K)$ and $0\leq k\leq n$ the set of all maps $h\in C(M\times K,\uin^k)$ with each $h|(\{z\}\times K)$, $z\in M$, being $(n-k,\gamma)$-discrete $($as a map from $K$ into $\uin^k$$)$ 
is dense in $C(M\times K,\uin^k)$ with respect to the uniform convergence topology.
\end{lem}

\begin{proof}
Suppose first that $K$ is metrizable and let $p_M\colon M\times K\to M$ and $p_K\colon M\times K\to K$ be the projections.
Then, by Hurewicz's theorem \cite{kk:68}, there exists a $0$-dimensional map 
$h^*\colon K\to\uin^n$. Consequently, $g^*=h^*\circ p_K$ is a map from $M\times K$ into $\uin^n$ such that $p_M\times g^*\colon M\times K\to M\times\uin^n$ is also $0$-dimensional.
According to Levin's \cite{ml:96} and Sternfeld's \cite{ys:95} results, the existence of such a map $g^*$ implies that the set ${\mathcal M}_n$ of all maps $g\in C(M\times K,\uin^n)$
with $\dim(p_M\times g)\leq 0$ is dense in $C(M\times K,\uin^n)$ with respect to the uniform convergence topology. If $q\colon\uin^n\to\uin^k$ is the projection generated by the first $k$ coordinates, then the map $g\rightarrow q\circ g$ is a continuous surjection from $C(M\times K,\uin^n)$ onto $C(M\times K,\uin^k)$ (both equipped with the uniform convergence topology), so
${\mathcal M}_k=\{q\circ g:g\in{\mathcal M}_n\}$ is dense in 
$C(M\times K,\uin^k)$. Moreover, since $\dim q=n-k$ and each $p_M\times g$, $g\in{\mathcal M}_n$, is $0$-dimensional, $\dim (p_M\times h)\leq n-k$ for any $h\in{\mathcal M}_k$ (the last conclusion is implied by the Hurewicz theorem on closed maps which lower dimension \cite{wh:27}). Therefore, $h_z=h|(\{z\}\times K)$ is an $(n-k)$-dimensional map for every $z\in M$ and
$h\in{\mathcal M}_k$.
Let show that any such $h_z$ is $(n-k,\gamma)$-discrete. Indeed, for fixed $y\in h_z(K)$ we have $\dim h_z^{-1}(y)\leq n-k$. So, there exists $\nu\in cov(h_z^{-1}(y))$ refining $\gamma$  such that $ord(\nu)\leq n-k+1$. Applying Lemma 2.1, we obtain an $(n-k,\gamma)$-discrete set $W_y$ in $K$ which contains $h_z^{-1}(y)$. Finally, choose a neighborhood $V_y$ of $y$ in $\uin^k$ such that $h_z^{-1}(V_y)\subset W_y$ and observe that $h_z^{-1}(V_y)$ is 
$(n-k,\gamma)$-discrete.   

Suppose now $K$ is not metrizable and fix $\delta >0$ and $h_0\in C(M\times K,\uin^k)$. We are going to find $h\in C(M\times K,\uin^k)$ satisfying the requirement of the lemma and such that $h$ is $\delta$-close to $h_0$. To this end,
represent $K$ as the limit space of a $\sigma$-complete inverse system $\mathcal S=\{K_{\beta},\pi_{\beta}^{\beta +1}:\beta\in B\}$ such that each $K_{\beta}$ is a metrizable compactum with $\dim K_{\beta}\leq n$. 
Applying standard inverse spectra arguments (see \cite{book}), we can find $\theta\in B$, $\gamma _1\in cov(K_{\theta})$ and $h_{\theta}\in C(M\times K_{\theta},\uin^k)$ such that $h_{\theta}\circ (id_M\times\pi_{\theta})=h_0$ and $\pi_{\theta}^{-1}(\gamma _1)$ refines $\gamma$, where $\pi_{\theta}\colon K\to K_{\theta}$ denotes the $\theta$th limit projection. Then, by virtue of the previous case, there exists a map 
$h_1\in C(M\times K_{\theta},\uin^k)$ which is $\delta$-close to $h_{\theta}$ and 
$h_1|(\{z\}\times K_{\theta})$ is $(n-k,\gamma_1)$-discrete. It follows from our construction 
that $h=h_1\circ (id_M\times\pi_{\theta})$ is $\delta$-close to $h_0$ and 
$h|(\{z\}\times K)$ is $(n-k,\gamma)$-discrete.
\end{proof}

Recall that a closed subset $F$ of the metrizable apace $M$ is said to be a $Z$-set in $M$,   if the set $C(Q,M\backslash F)$ is dense in $C(Q,M)$ with respect to the uniform convergence topology, where $Q$ denotes the Hilbert cube. If, in the above definition, $Q$ is replaced by $\uin^m$, $m\in\N\cup\{0\}$, we say that
$F$ is a $Z_m$-set in $M$.

\begin{lem}
Let $\alpha\colon X\to (0,\infty)$ be a positive continuous function and $g_0\in C(X,\uin^k)$.
Then $\psi_{(n-k,\omega)}(y)\cap\overline{B}(g_0,\alpha)$ is a $Z$-set in $\overline{B}(g_0,\alpha)$ for every $y\in Y$, where $\overline{B}(g_0,\alpha)$ is considered as a subspace of $C(X,\uin^k)$ with the uniform convergence topology. 
\end{lem}

\begin{proof}
In this proof all function spaces are equipped with the uniform convergence topology.
Since, by Lemma 2.6, $\psi_{(n-k,\omega)}$ has a closed graph, each $\psi_{(n-k,\omega)}(y)$
is closed in $C(X,\uin^k)$. Hence, $\psi_{(n-k,\omega)}(y)\cap\overline{B}(g_0,\alpha)$ is closed in $\overline{B}(g_0,\alpha)$. We need to show that, for fixed $y\in Y$, $\delta>0$ and a map $u\colon Q\to \overline{B}(g_0,\alpha)$ there exists a map
$v\colon Q\to\overline{B}(g_0,\alpha)\backslash\psi_{(n-k,\omega)}(y)$ which is $\delta$-close to $u$ with respect to the uniform metric $D_k$. To this end, observe first that $u$ generates $h\in C(Q\times X,\uin^k)$, $h(z,x)=u(z)(x)$, such that
$d_k(h(z,x),g_0(x))\leq\alpha (x)$ for any $(z,x)\in Q\times X$. Since $f^{-1}(y)$ is compact, we can find $\lambda\in (0,1)$ such that $\lambda\sup\{\alpha (x):x\in f^{-1}(y)\}<\displaystyle\frac{\delta}{2}$. Now, define $h_1\in C(Q\times f^{-1}(y),\uin^k)$ by $h_1(z,x)=(1-\lambda)h(z,x)+\lambda g_0(x)$. Then, for every $(z,x)\in Q\times f^{-1}(y)$, we have \\

\smallskip\noindent   
(3) \hbox{}~~~~~~$d_k(h_1(z,x),g_0(x))\leq (1-\lambda)\alpha (x)<\alpha (x)$ \\

\smallskip\noindent
and

\smallskip\noindent
(4) \hbox{}~~~~~~$d_k(h_1(z,x),h(z,x))\leq\lambda\alpha (x)<\displaystyle\frac{\delta}{2}$.

\smallskip\noindent
Let $\displaystyle q<\min\{r,\frac{\delta}{2}\}$, where $r$ is the positive number
$\inf\{\alpha (x)-d_k(h_1(z,x),g_0(x)):(z,x)\in Q\times f^{-1}(y)\}$. 
Since $\dim f^{-1}(y)\leq n$, by Lemma 2.7 (applied to the product $Q\times f^{-1}(y)$),   
there is a map $h_2\in C(Q\times f^{-1}(y),\uin^k)$ such that $d_k(h_2(z,x),h_1(z,x))<q$ and $h_2|(\{z\}\times f^{-1}(y))$ is an $(n-k,\omega)$-discrete map for each $(z,x)\in Q\times f^{-1}(y)$. Then, by $(3)$ and $(4)$, for all $(z,x)\in Q\times f^{-1}(y)$ we have \\

\smallskip\noindent
(5) \hbox{}~~~~~~$d_k(h_2(z,x),h(z,x))<\delta$ and $d_k(h_2(z,x),g_0(x))<\alpha (x)$. \\

\smallskip\noindent
Because both $Q$ and $f^{-1}(y)$ are compact, $u_2(z)(x)=h_2(z,x)$ defines 
the map $u_2\colon Q\to C(f^{-1}(y),\uin^k)$. The required map $v$ will be obtained as a lifting of $u_2$.
The restriction map $\pi\colon\overline{B}(g_0,\alpha)\to C(f^{-1}(y),\uin^k)$, $\pi (g)=g|f^{-1}(y)$, is obviously continuous (with respect to the uniform convergence topology).

\medskip
{\em Claim.} {\em $\pi\colon\overline{B}(g_0,\alpha)\to\pi (\overline{B}(g_0,\alpha))$ is an open map}. 

\medskip 
It's enough to show that \\

\noindent
(6) \hbox{}~~~~~~$\pi (\overline{B}(g_0,\alpha)\cap B_{\epsilon}(g))=\pi (\overline{B}(g_0,\alpha))\cap B_{\epsilon}(\pi(g))$ \\

\noindent
for every 
$g\in \overline{B}(g_0,\alpha)$ and $\epsilon>0$, where $B_{\epsilon}(g)$ and 
$B_{\epsilon}(\pi(g))$ are open balls, respectively, in 
$C(X,\uin^k)$ and $C(f^{-1}(y),\uin^k)$, both with the uniform metric generated by $d_k$. Let $p\in\pi (\overline{B}(g_0,\alpha))\cap B_{\epsilon}(\pi(g))$. Then $d_k(p(x),g_0(x))\leq\alpha(x)$ and $d_k(p(x),g(x))<\eta<\epsilon$ for every $x\in f^{-1}(y)$ and some positive number $\eta$. 
Define the closed and convex-valued map $\Phi\colon X\to {\mathcal F}_c(\uin^k)$ by
$\Phi (x)=p(x)$ if $x\in f^{-1}(y)$ and 
$\Phi (x)=\overline{B_{\alpha (x)}(g_0(x))\cap B_{\eta}(g(x))}$ if $x\not\in f^{-1}(y)$ (recall that $B_{\alpha (x)}(g_0(x))$ and $B_{\eta (x)}(g(x))$ are open balls in $\uin^k$). Since $g\in \overline{B}(g_0,\alpha)$, $\Phi(x)\neq\emptyset$ for every $x\in X$. Moreover, since $\alpha$, $g$ and $g_0$ are continuous, $\Phi$ is lower semi-contnuous. Therefore, by Michael's convex-valued selection theorem, $\Phi$ admits a selection $g_1\in C(X,\uin^k)$. Then $\pi(g_1)=p$ and $g_1\in\overline{B}(g_0,\alpha)\cap B_{\epsilon}(g)$. So, 
$\pi (\overline{B}(g_0,\alpha))\cap B_{\epsilon}(\pi(g))\subset\pi (\overline{B}(g_0,\alpha)\cap B_{\epsilon}(g))$ and, because the converse inclusion is trivial, we are done.     

\smallskip
Before completing the final step of our proof, note that 
$u_2(z)\in\pi (\overline{B}(g_0,\alpha))$ for every $z\in Q$ (indeed, consider the lower semi-continuous set-valued map $\phi\colon X\to{\mathcal F}_c(\uin^k)$, 
$\phi(x)=u_2(z)(x)$ for $x\in f^{-1}(y)$ and $\phi(x)=\overline{B_{\alpha (x)}(g_0(x))}$ for $x\not\in f^{-1}(y)$, and take any continuous selection $g$ of $\phi$). Now, we are going to lift the map $u_2$ to a map $v\colon Q\to\overline{B}(g_0,\alpha)$ such that $v$ is $\delta$-close to $u$. To this end,  
define $\theta\colon Q\to{\mathcal F}_c(C(X,\uin^k))$ by 
$\theta(z)=\overline{\pi^{-1}(u_2(z))\cap B_{\delta}(u(z))}$. The first inequality in $(5)$ iplies that $u_2(z)\in B_{\delta}(\pi(u(z)))$ for every $z\in Q$. Since each $u_2(z)$ belongs to $\pi(\overline{B}(g_0,\alpha))$, by virtue of $(6)$, $\theta(z)\neq\emptyset$, $z\in Q$.
On the other hand, since $\pi$ is open, by \cite[Example $1.1^*$ and Proposition 2.5]{em:56}, $\theta$ is lower semi-continuous. Obviously, every image $\theta(z)$ is convex and closed in $C(X,\uin^k)$ which is, in its own turn, closed and convex in the Banach space of all bounded and continuous functions from $X$ into $\R^k$. Therefore, using again the Michael selection theorem \cite[Theorem 3.2"]{em:56}, we can find a continuous selection $v\colon Q\to C(X,\uin^k)$ for $\theta$. Then $v$ maps $Q$ into $\overline{B}(g_0,\alpha)$ and $v$ is $\delta$-close to $u$. Moreover, for any $z\in Q$ we have $\pi(v(z))=u_2(z)$ and $u_2(z)$, being the restriction $h_2|(\{z\}\times f^{-1}(y))$, is $(n-k,\omega)$-discrete. Hence, 
$v(z)\in C_{(n-k,\omega)}(X|f^{-1}(y),\uin^k)$, $z\in Q$, i.e.
$v(Q)\subset\overline{B}(g_0,\alpha)\backslash\psi_{(n-k,\omega)}(y)$.       
\end{proof}

\begin{lem}
If $Y$ is a $C$-space, then $C_{(n-k,\omega)}(X,\uin^k)$ is dense in $C(X,\uin^k)$ with respect to the sourse limitation topology.
\end{lem}

\begin{proof}
It suffices to show that, for fixed $g_0\in C(X,\uin^k)$ and a positive continuous function $\alpha\colon X\to (0,\infty)$, there exists $g\in \overline{B}(g_0,\alpha)\cap C_{(n-k,\omega)}(X,\uin^k)$. We equipp $C(X,\uin^k)$ with the uniform convergence topology and 
consider
the constant (and hence, lower semi-continuous) map $\phi\colon Y\to{\mathcal F}_c(C(X,\uin^k))$, 
$\phi(y)=\overline{B}(g_0,\alpha)$. According to Lemma 2.8, 
$\overline{B}(g_0,\alpha)\cap\psi_{(n-k,\omega)}(y)$ is a $Z$-set in $\overline{B}(g_0,\alpha)$ for every $y\in Y$. So, we have a lower semi-continuous map $\phi\colon Y\to{\mathcal F}_c(E)$ and a map $\psi_{(n-k,\omega)}\colon Y\to 2^E$ such that $\psi_{(n-k,\omega)}$ has a closed graph (see Lemma 2.6) and $\phi(y)\cap\psi_{(n-k,\omega)}(y)$ is a $Z$-set in $\phi(y)$ for each $y\in Y$, where $E$ is the Banach space of all bounded continuous maps from $Y$ into $\R^k$. Therefore, we can apply \cite[Theorem 1.1]{gv:99} to obtain a  continuous map $h\colon Y\to E$ with $h(y)\in\phi(y)\backslash\psi_{(n-k,\omega)}(y)$ for every $y\in Y$ (Theorem 1.1 from \cite{gv:99} was proved under the assumption that $\psi_{(n-k,\omega)}$ has non-empty values, but the proof given in \cite{gv:99} works without this restriction). Observe that $h$ is a map from $Y$ into $\overline{B}(g_0,\alpha)$ such that $h(y)\not\in\psi_{(n-k,\omega)}(y)$ for every $y\in Y$, i.e. 
$h(y)\in \overline{B}(g_0,\alpha)\cap C_{(n-k,\omega)}(X|f^{-1}(y),\uin^k)$, $y\in Y$. Then
$g(x)=h(f(x))(x)$, $x\in X$, defines a map $g\in \overline{B}(g_0,\alpha)$ such that 
$g\in C_{(n-k,\omega)}(X|f^{-1}(y),\uin^k)$ for every $y\in Y$. Hence, by virtue of Corollary 2.4, $g\in C_{(n-k,\omega)}(X,\uin^k)$.  
\end{proof} 

%%%%%%%%%%%%%%%%%%%%%%%%%%%%%%%%%

%%%%%%%%%%%%%%%%%%%%%%%%%%%%%%
%%%%%%%%%%%%%%%%%%%%%%%%%%%%%%

\section{Proof ot Theorem 1.3}

The following proposition proves Theorem 1.3 in the special case when $f$ is $\sigma$-perfect.

\begin{pro}
Let $f\colon X\to Y$ be a $\sigma$-perfect map of metrizable spaces with $\dim f\leq n$ and $Y$ being a $C$-space. Then the set of all maps $g\colon X\to \uin^n$ such that $\dim (f\times g)=0$ is dense and $G_{\delta}$ in $C(X,\uin^n)$ with respect to the source limitation topology.
\end{pro}

\begin{proof}
All function spaces in this proof are considered with the source limitation topology.
Let $X$ be the union of the closed sets $X_i$, $i=1,2,..$, such that each restriction $f_i=f|X_i$ is perfect and $Y_i=f(X_i)$ is closed in $Y$. Fix a sequence $\{\omega_q\}$ of open covers of $X$ with $mesh(\omega_q)<q^{-1}$. Every $Y_i$ is a $C$-space (as a closed set in $Y$), so we can apply Lemma 2.9 to the maps $f_i\colon X_i\to Y_i$ and conclude that  
$\displaystyle H_i=\bigcap_{q=1}^{\infty}C_{(0,\omega_q)}(X_i,\uin^n)$ is dense and $G_{\delta}$ in $C(X_i,\uin^n)$, $i\in\N$. Here, $C_{(0,\omega_q)}(X_i,\uin^n)$ consists of all
$h\in C(X_i,\uin^n)$ such that $f_i\times h$ is $(0,\omega_q)$-discrete.
Since all restriction maps $p_i\colon C(X,\uin^n)\to C(X_i,\uin^n)$, $p_i(g)=g|X_i$, are continuous, open and surjective, the sets $C_i=p_i^{-1}(H_i)$ are dense and $G_{\delta}$ in $C(X,\uin^n)$, so is the intersection $\bigcap_{i=1}^{\infty}C_i$. It only remains to observe that $g\in\bigcap_{i=1}^{\infty}C_i$ if and only if $\dim(f_i\times g_i)=0$ for every $i$, where $g_i=g|X_i$. Hence, by the countable sum theorem, 
$g\in\bigcap_{i=1}^{\infty}C_i$ if and only if $\dim(f\times g)=0$.
\end{proof} 

We continue now with the proof of the first part of Theorem 1.3.
Suppose $f\colon X\to Y$ is a closed $n$-dimensional surjection with both $X$ and $Y$ metrizable and $Y$ a $C$-space. By Vainstein lemma \cite{re:89}, the boundary 
$Fr f^{-1}(y)$ of every $f^{-1}(y)$ is compact. Defining $F(y)$ to be $Fr f^{-1}(y)$ if $Fr f^{-1}(y)\ne\emptyset$, and an arbitrary point from $f^{-1}(y)$ otherwise, we obtain the set $X_0=\cup\{F(y):y\in Y\}$ such that $X_0\subset X$ is closed and the restriction $f|X_0\colon X\to Y$ is a perfect surjection. Moreover, each $f^{-1}(y)\backslash X_0$ is open in $X$, so $\dim(X\backslash X_0)\leq n$. Represent $X\backslash X_0$ as the union of countably many closed sets
$X_i\subset X$ and for each $i=0,1,2,..$ let $p_i\colon C(X,\uin^n)\to C(X_i,\uin^n)$ be the restriction map. By Proposition 3.1, the set $C_0$ consisting of all $g\in C(X,\uin^n)$ with $(f\times g)|X_0$ being $0$-dimensional is dense and $G_{\delta}$ in $C(X,\uin^n)$ with respect to the source limitation topology. Consequently, $C_0$ is uniformly dense in $C(X,\uin^n)$. On the other hand, 
since $\dim X_i\leq n$ for every $i=1,2,..$, the set $H_i\subset C(X_i,\uin^n)$ of all uniformly $0$-dimensional maps is dense and $G_{\delta}$ in $C(X_i,\uin^n)$ with respect to the uniform convergence topology \cite{mk:52} (recall that a map $h\colon X_i\to\uin^n$ is uniformly $0$-dimensional if for every $\epsilon >0$ there exists $\eta>0$ such that, if $T\subset\uin^n$ and $diam(T)\leq\eta$, then $h^{-1}(T)$ is covered by a disjoint open family in $X_i$ consisting of sets with diameter $\leq\epsilon$). Because $p_i$ are open and continuous surjections when $C(X,\uin^n)$ and $C(X_i,\uin^n)$ carry the uniform convergence topology, 
all $C_i=p_i^{-1}(H_i)$, $i=1,2,..$, are uniformly dense and $G_{\delta}$ in $C(X,\uin^n)$.  Therefore, $C_{\infty}=\bigcap_{i=0}^{\infty}C_i$ is $G_{\delta}$ in $C(X,\uin^n)$ with respect to the source limitation topology. Moreover, $f\times g$ is $0$-dimensional for every $g\in C_{\infty}$. It remains to show that $C_{\infty}$ is uniformly dense in $C(X,\uin^n)$.
For every $g\in C_0$ let
$H(g)=\{h\in C(X,\uin^n): h|X_0=g|X_0\}$. Obviously, $C_0=\cup\{H(g):g\in C_0\}$ and each $H(g)$ is uniformly closed in $C(X,\uin^n)$. So, $C_{\infty}$ is the union of the sets 
$A(g)=\bigcap_{i=1}^{\infty}C_i\cap H(g)$, $g\in C_0$. For fixed $g\in C_0$ and $i=1,2,..$, let $p_i(g)=p_i|H(g)$. Using that $X_0$ and $X_i$ are closed disjoint subsets of $X$, one can show that every $p_i(g)\colon H(g)\to C(X_i,\uin^n)$ is an uniformly continuous and open surjection. Hence, $H(g)\cap C_i$ is dense and $G_{\delta}$ in $H(g)$ with respect to the uniform convergence topology as the preimage of $H_i$ under $p_i(g)$. Therefore, $A(g)$ is uniformly dense in $H(g)$ (recall that $H(g)$ is uniformly closed in $C(X,\uin^n)$, so it has Baire property). We finally observe that the uniform density of $C_0$ in $C(X,\uin^n)$ and the uniform density of $A(g)$ in $H(g)$ for each $g\in C_0$ yield the uniform density of $C_{\infty}$ in $C(X,\uin^n)$.

%%%%%%%%%%%%%%%%%%%%%%%%%%%%%%%%%

%%%%%%%%%%%%%%%%%%%%%%%%%%%%%%
%%%%%%%%%%%%%%%%%%%%%%%%%%%%%%

\section{Proof ot Theorem 1.4}

It suffices to prove this theorem for closed maps, so we suppose that $f\colon X\to Y$ is a closed surjection. If $A_{n-1}$ is constructed, then for $k<n-1$, we can find an $F_{\sigma}$-subset $A_k\subset A_{n-1}$ with $\dim A_k\leq k$ and $\dim(A_{n-1}\backslash A_k)\leq n-k-2$ (this can be accomplished by induction, the first step is to represent $A_{n-1}$ as the union of $0$-dimensional $G_{\delta}$-subsets $B_j$, $j=1,2,..n$ and to denote $A_{n-2}=\cup_{j=1}^{j=n-1}B_j$). Therefore, we need only to construct $A_{n-1}$. To this end, we first establish the following analogue of Sternfeld's \cite[Lemma 1]{ys:95} which was proved for compact metrizable spaces.

\begin{lem}
Let $M$ be metrizable and $K$ a compact metric space with $\dim K\leq n$. Then there exists a $F_{\sigma}$ subset $B\subset M\times K$ such that $\dim B\leq n-1$ and $\pi_M|(M\times K)\backslash B$ is $0$-dimensional, where $\pi_M\colon M\times K\to M$ is the projection.
\end{lem} 

\begin{proof}
As in \cite{ys:95}, the proof can be reduced to the case $n=1$ and $K=[0,1]$. So, we are going to show the existence of a $0$-dimensional $F_{\sigma}$-subset $B$ of $M\times\uin$ such that each set $(\{y\}\times\uin)\backslash B$, $y\in M$, is $0$-dimensional and that will complete the proof. Let 
$h\colon Z\to M$ be a perfect surjection with $Z$ being a $0$-dimensional metrizable space. Then, by \cite[Proposition 9.1]{bp:98}, there exists a map $g\colon Z\to Q$ such that $h\times g\colon Z\to M\times Q$ is a closed embedding. Next, let $\Delta$ be the Cantor set and take a surjection $p\colon\Delta\to Q$ admitting an averaging operator between the function spaces $C(\Delta)$ and $C(Q)$ \cite{ap:68} (such maps are called Milyutin maps). According to \cite{sd:73}, there exists a lower semi-continuous compact-valued map $\phi\colon Q\to 2^{\Delta}$ with $\phi(y)\subset p^{-1}(y)$ for every $y\in Q$. 
We can apply Michael's $0$-dimensional selection theorem \cite{em:57} to obtain a continuous selection $q$ for the map $\phi\circ g$.  Obviously $h\times q\colon Z\to M\times\Delta$ is a closed embedding, so $Z_0=(h\times q)(Z)$ is a $0$-dimensional closed subset of $M\times\Delta$. Finally, considering $\Delta$ as a subset of $\uin$, let $Z_r=\{(h(z),q(z)+r):z\in Z\}\subset M\times\uin$ for every rational $r\in\uin$, where addition $q(z)+r$ is taken in $\R$ $mod1$.
Then each $Z_r$ is a closed subset of $M\times\uin$ homeomorphic to $Z$, so $B=\cup\{Z_r:\hbox{$r$ is rational}\}$ is $0$-dimensional and $F_{\sigma}$ in $M\times\uin$. Moreover,  $(\{y\}\times\uin)\backslash B$ is also $0$-dimensional for every $y\in M$.   
\end{proof}

Let continue the proof of Theorem 1.4. As in the proof of Theorem 1.3, there are closed subsets $X_i\subset X$, $i=0,1,2,..$, such that $f|X_0$ is a perfect map onto $Y$, each $X_i$, $i\geq 1$, is at most $n$-dimensional and $X\backslash X_0=\bigcup_{i=1}^{\infty}X_i$. 
For every $i\geq 1$ we choose an $(n-1)$-dimensional $F_{\sigma}$-set $H_i\subset X_i$ with $\dim (X_i\backslash H_i)\leq 0$. A similar type subset of $X_0$ can also be found. Indeed, let $f_0=f|X_0$ and take $g\colon X_0\to\uin^n$ such that $f_0\times g\colon X_0\to Y\times\uin^n$ is $0$-dimensional (see Theorem 1.3). By Lemma 4.1, there exists an $F_{\sigma}$-set $B\subset Y\times\uin^n$ with $\dim B\leq n-1$ and each $(\{y\}\times\uin^n)\backslash B$, $y\in Y$, being $0$-dimensional. Then $H_0=(f_0\times g)^{-1}(B)$ is $F_{\sigma}$ in $X_0$.
Since $f_0\times g$ is perfect, by the generalized Hurewicz theorem on closed maps lowering dimension \cite{es:62}, we have $\dim H_0\leq n-1$ and $\dim(f_0^{-1}(y)\backslash H_0)\leq 0$ for every $y\in Y$. Finally, set $A_{n-1}=\bigcup_{i=0}^{\infty}H_i$. Obviously, $\dim A_{n-1}\leq n-1$. On the other hand, each $f^{-1}(y)\backslash A_{n-1}$, $y\in Y$, is the union of its closed sets $F_i(y)=f^{-1}(y)\cap X_i\backslash A_{n-1}$, $i\geq 0$. But $F_0(y)=f_0^{-1}(y)\backslash H_0$ and $F_i(y)\subset f^{-1}(y)\cap(X_i\backslash H_i)$ for $i\geq 1$, so all $F_i(y)$ are $0$-dimensional. Consequently, $\dim(f^{-1}(y)\backslash A_{n-1})\leq 0$ for every $y\in Y$, i.e. the restriction $f|(X\backslash A_{n-1})$ is $0$-dimensional.

%%%%%%%%%%%%%%%%%%%%%%%%%%%%%%%%%

%%%%%%%%%%%%%%%%%%%%%%%%%%%%%%
%%%%%%%%%%%%%%%%%%%%%%%%%%%%%%

\section{Some applications.}

Our first application deals with extensional dimension introduced by Dranishnikov \cite{d:95} (see also \cite{ch:00} and \cite{dd:95}). Let $K$ be a $CW$-complex and $X$ a normal space.  We say that the extensional dimension of $X$ doesn't exceed $K$, notaion $\e X\leq K$, if every map $h\colon A\to K$, where $A\subset X$ is closed, can be extended to a map from $X$ into $K$ provided there exist a neighborhood $U$ of $A$ in $X$ and a map $g\colon U\to K$ extending $h$. Obviously, if $K$ is an absolute neighborhood extensor for $X$, then $\e X\leq K$ iff $K$ is an absolute extensor for $X$. In this notation, $\dim X\leq n$ is equivalent to $\e X\leq\sphere^n$. We also write $\e X\leq\e Y$ if $\e Y\leq K$ implies $\e X\leq K$ for any $CW$-complex $K$. 

Dranishnikov and Uspenskij \cite{du:97} provided a generalization of the Hurewicz theorem on dimension lowering maps: if $f\colon X\to Y$ is an $n$-dimensional surjection between compact finite-dimensional spaces, then $\e X\leq\e(Y\times\uin^n)$; moreover, this statement holds for any compact spaces (not necessary finite-dimensional) when $n=0$. 
We can improve this result as follows (see also \cite{cv:99} and \cite{drs:98} for another extension dimensional variants of Hurewicz's theorem):

\begin{thm}
If $f\colon X\to Y$ is a perfect $n$-dimensional surjection such that $Y$ is a paracompact $C$-space, then $\e X\leq\e(Y\times\uin^n)$.
\end{thm}

Theorem 5.1 follows from Theorem 2.2 and next proposition which can be extracted from the Dranishnikov and Uspenskij proof of their \cite[Lemma 2.1 and Theorem 1.4]{du:97}.

\begin{pro}
Let $K$ be a $CW$-complex and $X$ paracompact. If for any $\omega\in cov(X)$ there exist a paracompact space $Z_{\omega}$ with $\e Z_{\omega}\leq K$ and a perfect $(0,\omega)$-discrete map $g\colon X\to Z_{\omega}$, then $\e X\leq K$.
\end{pro}

\begin{cor}
Let $f\colon X\to Y$ be a $\sigma$-closed $n$-dimensional surjection between metrizable spaces with $Y$ being a $C$-space. Then $\e X\leq\e (Y\times\uin^n)$.
\end{cor}

\begin{proof}
Let $K$ be a $CW$-complex with $\e(Y\times\uin^n)\leq K$. It suffices to show that $\e X\leq K$. Since extension dimension satisfies the countable sum theorem, the proof of the last inequality is reduced to the case $f$ is closed. We can also assume that $K$ is an open subset
of a normed space because every
$CW$-complex is homotopy equivalent to such a set. 
Represent $X$ as the union of the closed sets $X_i\subset X$, $i\geq 0$, such that $f|X_0$ is a perfect map onto $Y$ and $\dim X_i\leq n$ for each $i\geq 1$ (see the proof of Theorem 1.3). Then, by Theorem 5.1, $\e X_0\leq K$. On the other hand, $\e (Y\times\uin^n)\leq K$ implies that $\e\uin^n\leq K$, in particular, every map from $\sphere^{n-1}$ into $K$ is extendable to a map from $\uin^n$ into $K$. In other words, $K$ is $C^{n-1}$ and, as an open subset of a normed space, $K$ is also $LC^{n-1}$. It is well known that $LC^{n-1}$ and $C^{n-1}$ metrizable spaces are precisely the absolute extensors for $n$-dimensional metrizable spaces. Hence, $\e X_i\leq K$ for every $i\geq 1$. Finally, by the countable sum theorem for extensional dimension, we have $\e X\leq K$.
\end{proof}
 
Another application is a parametric version of the Bogatyi decomposition theorem of $n$-dimensional metrizable spaces \cite{sb:99}: For every metrizable $n$-dimensional space $M$ there exist countably many $0$-dimensional $G_{\delta}$-subsets $M_k\subset M$ such that 
$M=\bigcup_{i=1}^{i=n+1}M_{k(i)}$ for all pairwise distinct $k(1),..,k(n+1)$ in $\N$.

\begin{pro}
Let $f\colon X\to Y$ be a closed $n$-dimensional surjection between metrizable spaces with $Y$ a $C$-space. Then there exists a sequence $\{A_k\}$ of $G_{\delta}$-subsets of $X$ such that
every restriction $f|A_k$ is $0$-dimensional and for any $P\subset\N$ of cardinality $n+1$ we have $X=\bigcup_{k\in P}A_k$. 
\end{pro}

\begin{proof}
Take closed sets $X_i\subset X$, $i\geq 0$, and a map $g\colon X\to\uin^n$ such that $f|X_0$ is perfect, $X\backslash X_0=\bigcup_{i\geq 1}X_i$, $\dim (f\times g)=0$ and each $g|X_i$, $i\geq 1$, is uniformly $0$-dimensional (see the proof of Theorem 1.3). According to the Bogatyi theorem, there exists a sequence of $0$-dimensional $G_{\delta}$-subsets $B_k\subset\uin^n$ such that $\uin^n$ is the union of any $n+1$ elements of this sequence. Let $A_k=(f\times g)^{-1}(Y\times B_k)$, $k\in\N$. The only non-trivial condition we need to check is that each restriction $f|A_k$ is $0$-dimensional, i.e. $\dim f^{-1}(y)\cap A_k\leq 0$ for all $y\in Y$ and $k\geq 1$. For fixed $y$ and $k$ we have $f^{-1}(y)\cap A_k=\bigcup_{i\geq 0}g_i^{-1}(B_k)$, where $g_i$ denotes the restriction $g|(f^{-1}(y)\cap X_i)$.
Since every $g_i^{-1}(B_k)$ is closed in 
$f^{-1}(y)\cap A_k$, it suffices to show that the sets $g_i^{-1}(B_k)$, $i\geq 0$, are $0$-dimensional. For $i=0$ this follows from the Hurewicz lowering dimension theorem \cite{wh:27} because $g_0$ is a perfect $0$-dimensional map. For $i\geq 1$ we use that $g|X_i$ is uniformly $0$-dimensional and the preimage of any $0$-dimensional set under uniformly $0$-dimensional map is again $0$-dimensional.  
\end{proof}

A map $f\colon X\to Y$ is said to be of countable functional weight \cite{bp:98} (notation $W(f)\leq\aleph_0$) if there exists a map $h\colon X\to Q$, $Q$ is the Hilbert cube, such that 
$f\times h\colon X\to Y\times Q$ is an embedding. In \cite{bp:98} Pasynkov has shown that his  results from \cite{bp:96} remain valid for maps $f\colon X\to Y$ between finite-dimensional completely regular spaces $X$ and $Y$ such that $W(f)\leq\aleph_0$ and both $f$ and its \v{C}ech-Stone extension have the same dimension (the last condition holds, for example, if $X$ is normal, $Y$ paracompact and $f$ closed). We are going to show that Theorem 2.2 implies a similar result with $Y$ being a $C$-space.

\begin{thm}
Let $f\colon X\to Y$ be a $\sigma$-closed $n$-dimensional surjection of countable functional weight such that $X$ is normal and $Y$ a paracompact $C$-space. Then the set $\mathcal G$ of all maps $g\in C(X,\uin^n)$ with
$\dim(f\times g)=0$ is uniformly dense in $C(X,\uin^n)$. If, in addition, $X$ is paracompact and $f$ is $\sigma$-perfect, then $\mathcal G$ is dense and $G_{\delta}$ in $C(X,\uin^n)$ with respect to the source limitation topology.
\end{thm}    

\begin{proof}
Since $W(f)\leq\aleph_0$, there exists a map $h\colon X\to Q$ such that $f\times h$ is an embedding. For every $k\in\N$ let $\gamma_k$ be an open cover of $Q$ of mesh $\leq k^{-1}$. 
Suppose $f$ is $\sigma$-closed and represent $X$ and $Y$ as the union of closed sets $X_i$ and $Y_i$, respectively, such that each $f_i=f|X_i$ is a closed map onto $Y_i$. Let
$Z_i=(\beta f_i)^{-1}(Y_i)$ and $\bar{f_i}=(\beta f_i)|Z_i$, $i\in\N$, where $\beta f_i$ denotes the \v{C}ech-Stone extension of $f_i$.  
Because $X$ is normal, each $Z_i$ is a closed subset of $Z=(\beta f)^{-1}(Y)$ and 
 $\bar{f_i}\colon Z_i\to Y_i$ are perfect $n$-dimensional maps. Moreover, $Z$ is paracompact as the preimage of $Y$. We consider the extension $\bar{h}\colon Z\to Q$ of $h$ and the covers $\omega_k=\bar{h}^{-1}(\gamma_k)\in cov(Z)$. By Theorem 2.2 (applied for the maps $\bar{f_i}$), the sets $\mathcal H_{i,k}$ consisting of all $g\in C(Z,\uin^n)$ such that $f_i\times g$ is $(0,\omega_k)$-discrete, $i,k\in\N$, are open and dense in $C(Z,\uin^n)$
with respect to the source limitation topology, so is $\mathcal H=\bigcap_{i,k=1}^{\infty}\mathcal H_{i,k}$. Moreover, $\mathcal H$ is uniformly dense in $C(Z,\uin^n)$.
Therefore, the set $\mathcal G_0=\{g|X:g\in\mathcal H\}$ is uniformly dense in $C(X,\uin^n)$. Since $h$ is a homeomorphism on every fiber of $f$, $\mathcal G_0\subset\mathcal G$. Hence, $\mathcal G$ is also uniformly dense in $C(X,\uin^n)$.  

Let $f$ be $\sigma$-perfect and $X$ paracompact. Then substituting $\bar{f_i}=f_i$, $Z_i=X_i$ and $Z=X$ in the previous proof, we obtain that $\mathcal G$ coincides with $\mathcal H$.
\end{proof}

\begin{cor}
Let $f\colon X\to Y$ be a $\sigma$-closed $n$-dimensional surjection having second countable fibres. If $X$ is metrizable and $Y$ a paracompact $C$-space, then the set of all $g\in C(X,\uin^n)$ with $\dim (f\times g)=0$ is uniformly dense in $C(X,\uin^n)$.
\end{cor}

\begin{proof} Since $f$ is of countable functional weight (see \cite[Proposition 9.1]{bp:98}), this corollary follows from Theorem 5.5.
\end{proof}

We finally formulate the following result, its proof is similar to that one of Theorem 2.2.  
 
\begin{thm}
Let $f\colon X\to Y$ be a perfect surjection of countable functional weight with $Y$ a paracompact $C$-space. Then all maps $g\in C(C,Q)$ such that $f\times g$ is an embedding form a dense and $G_{\delta}$ subset of $C(X,Q)$ with respect to the source limitation topology.
\end{thm}

\begin{cor}
Let $f\colon X\to Y$ be a perfect surjection between metrizable spaces. If $Y$ is a $C$-space, then the set of all $g\in C(X,Q)$ with $f\times g$ being embedding is dense and $G_{\delta}$ in $C(X,Q)$ with respect to the source limitation topology.
\end{cor}

\bigskip

\end{document}